\documentclass[12pt]{article}
\usepackage[english]{babel}
\usepackage{amssymb,amsmath,graphics,tikz}

\usepackage{hyperref,tipa}

\parskip\medskipamount
\newtheorem{theorem}{Theorem}[section]

\newtheorem{remark}[theorem]{Remark}

\newtheorem{lemma}[theorem]{Lemma}
\newtheorem{cor}[theorem]{Corollary}
\def\<{\langle}
\def\>{\rangle}
\newcommand{\proof}{\emph{Proof.~}}

\newcommand{\dd}{\mathsf{d}}

\newcommand{\cA}{\mathcal{A}}

\newcommand{\cF}{\mathcal{F}}

\def\qed{{\hfill\hphantom{.}\nobreak\hfill$\Box$}}

\newcommand{\id}{\mathrm{id}}

\newcommand{\A}{\mathbb{A}}
\newcommand{\R}{\mathbb{R}}

\title{Descent of affine buildings - I. Large minimal angles}
\author{ Bernhard M\"uhlherr \and Koen Struyve\thanks{The second author is supported by  the Fund for Scientific Research --
Flanders (FWO - Vlaanderen)} \and  Hendrik Van Maldeghem }
\begin{document}
\maketitle
\begin{abstract}
In this two-part paper we prove an existence result for affine buildings arising from exceptional algebraic reductive groups. Combined with earlier results on classical groups, this gives a complete and positive answer to the conjecture concerning the existence of affine buildings arising from such groups defined over a (skew) field with a complete valuation, as proposed by Jacques Tits.

This first part lays the foundations for our approach and deals with the `large minimal angle' case.
\end{abstract}
\section{Introduction}
A central problem in Bruhat-Tits theory is to show the existence of the affine building associated with a semi-simple algebraic group
over a field with valuation. Large parts of \cite{Bru-Tit:72} and \cite{Bru-Tit:84} are devoted to this problem and a fairly general solution
is provided. There are in principal two methods to achieve this goal. The first is to show directly that
a valuation of the defining division ring can be used to produce a root datum with valuation. In \cite{Bru-Tit:72} this method
was used to give a complete solution for all classical groups. We also mention~\cite{Wei:09} and~\cite{Wei:*}, where each case of relative rank at least two is handled with one exception. It appears that this approach becomes 
difficult, if one wants to handle certain forms of exceptional groups. This is because there are $k$-forms of certain
exceptional groups of relative rank 1 for which there are not yet convenient descriptions in terms of algebraic
parameter systems (e.g. pseudo-quadratic forms, Jordan algebras).
In \cite{Bru-Tit:72} and \cite{Bru-Tit:84} the existence problem for the exceptional groups
was treated by using Galois descent. One knows that there is an affine building associated to any split group and constructs
the desired affine building as fixed point structure of a Galois group acting on the affine building associated with the
group over the separable closure. This approach is very powerful in the sense that it provides the desired results for a wide
range of cases.
In fact, these results cover all the cases which are relevant for applications as for instance
the case of finite residue fields. 
However, the methods do not work in full generality. One has to require that there is a Galois extension of the ground field
over which the group splits and in which the valuation is tamely ramified. The method has  been refined and generalized
by G. Rousseau in \cite{rousseauthesis} and by G. Prasad  and J. Yu in \cite{PrasadYu}. 

The idea of using Galois descent is also at the base of our approach. Our approach is almost purely geometric.  
The main advantage of our new approach with respect
to the earlier results in this direction is the following. It is not true, that the fixed point set of a Galois group acting on 
the Bruhat-Tits building is an affine building. This phenomenon was already observed in \cite{Bru-Tit:72} and studied in some detail in \cite{rousseauthesis}.
The conditions in the Galois descent results mentioned above are designed to ensure that the fixed point set is very close to being a building.
An example for such a condition is that the order of the finite Galois group is prime to the characteristic of the residue field (see \cite{PrasadYu}).
We do not ask conditions on the group or the behavior of the valuation in the extension field. Instead, we ask for a geometric
condition using the Tits diagram of the algebraic group in question. This enables us to deal with cases where the fixed point set
is not a building in the strong sense, but a thickened version of the desired affine building.

The key notion for our approach is the {\sl minimal angle} of a Tits diagram for an algebraic group.
If this minimal angle  is not too small we are able  to 
prove the existence of the  affine buildings associated with an algebraic group having this Tits diagram.
There is the well understood special case where the minimal
angle of the Tits diagram is $\pi$. This happens precisely when the group is quasi-split and it was already 
observed in \cite{Bru-Tit:72} that the fixed point set of the Galois group acting on the building over the separable
closure is again an affine building -- in fact the affine building associated with the corresponding quasi-split form.
As already mentioned, the fixed point set of a Galois group need not to be an 'honest' affine building. Nevertheless
it is reasonable to expect that it is always a thickened version of an affine building. We show that this is indeed the case if
the minimal angle of the Tits diagram is strictly greater than $\pi/3$. In fact, we show that there is a canonical convex subset of  the  fixed point set,
which is an affine building. Since the formulation of our result needs some preparation we refer to Section \ref{mainresult} for
the precise statement.
 Unfortunately, we are not able to extend our  result to the case where the angle is not strictly greater  
than $\pi/3$. Nevertheless, in Part II of the present paper, further building on the results of the present Part I, the second author is able to prove the existence of the desired affine buildings under the assumption
that the minimal angle is \emph{at least} $\pi/3$. His construction somehow `extracts' the desired affine building out of the fixed point set.

As in the earlier approaches to the existence problem, also our methods do not cover all cases, simply because
we have no idea what to do if the minimal angle of the diagram is smaller than $\pi/3$. However, the final results
obtained in Part II of this paper can be used to cover all the cases for which there had been no
existence proof so far. Thus, the existence of an affine building associated to a semi-simple algebraic
group is established in full generality, if one takes the earlier results mentioned above into account. This will be explained
in detail in Part II.

Finally we would like to point out two favorable aspects of our geometric approach to the existence problem.

First of all, our method is almost purely geometric, and as such, we might as well include the nondiscrete case. This is exactly what we will do. Also, our construction is canonical to the extent that the group stabilizing the subbuilding at infinity will act in a natural way on the affine building that we construct for it. 

Secondly, our geometric point of view allows for a generalization to buildings which do no longer arise from pure algebraic groups, but also from groups of mixed type (terminology as in \cite{Tit:74}, where the latter are introduced and constructed). The standard example is the construction of the exceptional Moufang quadrangles of type $\mathsf{F}_4$, see \cite{Mue-Mal:98}. But also in some buildings that do arise from algebraic groups, the geometric construction works whereas the algebraic one fails. A standard example is the construction of the perfect Ree-Tits generalized octagons from buildings of type $\mathsf{F}_4$ over a perfect field of characteristic 2 (see also~\cite{Hit-Kra-Wei:10}). The corresponding procedure on the algebraic group does not exist, but viewed as a Chevalley group, it has an analogue known as the  `Frobenius twist'. The previous example is the only one with such a twist producing buildings of rank strictly greater than one (and note that, in this case, we are dealing with $\mathbb{R}$-buildings which are never simplicial). One reason why this geometric point of view works is the fact that the fixed point structure of any group acting on a building is a convex simplicial subcomplex which is needed to obtain a subbuilding.

\textbf{Acknowledgement.}  The authors would like to thank Richard Weiss for bringing this problem to our attention, and for valuable comments and discussions on the subject.

\section{Preliminaries}

In this section we recall some definitions and fix notation. Except for Subsections~\ref{titsdiagrams},~\ref{section:angle}, and~\ref{section:classification},  the content of this section is standard. 
We roughly follow the first three chapters of \cite{Abr-Bro:08} and the first section of~\cite{Par:00} to which we refer for further information.

\subsection{Simplicial complexes}

A \emph{simplicial complex} $S$ on a set $X$ is a set of finite subsets of $X$ such that for each subset $x \in S$ and $y\subset x$, we also have that $y \in S$. We also ask that each singleton of $X$ is in $S$. The elements of $X$ are called the \emph{vertices}, the elements of $S$ are called \emph{simplices}. We will always assume that the order of simplices is bounded. 

A \emph{maximal simplex} of a simplicial complex $S$ on $X$, is a simplex of $S$ not contained in a larger simplex. Two maximal simplices of the same order are called \emph{adjacent} if they share a simplex of order one less (which is called a \emph{panel}).

A simplicial complex is a \emph{chamber complex} if for each two maximal simplices $C$ and $D$ there is a sequence $(C_0 =C,C_1,\dots,C_i=D)$ of maximal simplices, such that every two subsequent maximal simplices are adjacent. In this case the maximal simplices are called \emph{chambers}. Note that this implies that all the chambers are of the same order. 

\subsection{Coxeter complexes}\label{section:coxeter}
A \emph{Coxeter matrix} is an $n \times n$-matrix $M$ such that $m_{ii} =1$ for $i \in \{1,\dots,n \}$ and $m_{ij} = m_{ji} \in \{2,3, \dots, \infty \}$ for $i, j \in \{1,\dots,n \}$ and $i \neq j$.

The \emph{Coxeter group} associated with this matrix $M$ is the group $W$ with generators $S=\{s_1, \dots ,s_n\}$ and relations $(s_i s_j)^{m_{ij}} = e$, with $e$ the identity element of $W$ (if $m_{ij}=\infty$ we do not put any relation). The pair $(W,S)$ is called the \emph{Coxeter system}. Note that the elements of $S$ are involutions. 

Usually one uses a \emph{Coxeter diagram} to represent the Coxeter matrix. This diagram consists of $n$ vertices, one for each generator in $S$. If for two different generators $s_i$ and $s_j$ one has that $m_{ij}=2$, then there is nothing drawn between the associated vertices; if $m_{ij}=3$, then one draws a single edge, if $m_{ij}=4$, a double edge. If $m_{ij}>4$ one uses an edge labeled with $m_{ij}$. The \emph{rank} of a Coxeter system or diagram is the size of the set $S$.

To a Coxeter system one can associate a chamber complex $\Sigma_W$, called the \emph{Coxeter complex}. The Coxeter group acts sharply transitive on the chambers of the Coxeter complex, and there exists a \emph{type} function from the simplices of the simplicial complex to the subsets of $S$ (or equivalently to sets of nodes of the Coxeter diagram). 


If the Coxeter group $W$ is finite, it can be realized as a finite reflection group acting on a $n$-dimensional Euclidean space. The conjugates of elements in $S$ form reflections of this real Euclidean space. The hyperplanes they fix are called \emph{walls}, the closed half-spaces they border \emph{roots}. The set of all walls defines cones (called \emph{sector-faces}) which correspond to the simplices of the Coxeter complex. The maximal cones are called \emph{sectors}, the one less than maximal \emph{sector-panels}. If we take the intersection with a sphere centered at the fixed point of the reflection group one obtains a geometric realization of the Coxeter complex on a $(n-1)$-dimensional sphere. For this reason one also calls this case \emph{spherical}. 

If the Coxeter group can be realized as an affine reflection group, realizing the Coxeter complex as a  triangulation of the Euclidean space, then one calls the Coxeter system \emph{affine}.

A Coxeter system is \emph{irreducible} if the associated Coxeter diagram is connected.

\subsection{Spherical Coxeter systems and diagram automorphisms}\label{section:diag}

In this subsection we collect some results on diagram automorphisms of spherical Coxeter systems. Hence, throughout this subsection
all Coxeter systems are assumed to be spherical.

Let $(W,S)$ be a finite Coxeter system and consider its geometric realization on the sphere. The natural metric on the sphere yields an \emph{angular distance} between
any two vertices of the Coxeter complex.

If $M$ is the Coxeter diagram (over $I$) associated with $(W,S)$ then we define the \emph{angular distance of $i \in I$ with respect
to $M$} as  the minimal distance between two distinct vertices of type $i$ on the sphere.

Suppose $M$ is of rank 1, which means that $I=\{ i \}$, then the angular distance of $i$ with respect to $M$
is $\pi$; if $M$ is of rank 2 representing the dihedral group of order 2m, then the angular distance at both indices is $2\pi/m$.
If $M = \mathsf{A}_3$ labeled by $\{ 1,2,3 \}$ in linear order, then the angular distances at $1$ and $3$ are $\arccos(-1/3)$, whereas the angular distance at $2$ is $\pi/2$.

The angular distance at a node of a spherical Coxeter diagram is one ingredient of the definition of a
 minimal angle of a Tits diagram.  A second is the consideration
of diagram automorphisms which we consider now.

The following well known facts can be found  in \cite{Muh:93}:
Let $\Gamma$ be a subgroup of the automorphism group of a spherical Coxeter diagram $M$ over a set $I$.
Then $\Gamma$ acts naturally on the corresponding Coxeter system $(W,S)$  and its associated Coxeter complex.
The centralizer $W_{\Gamma}$ of $\Gamma$ in $W$ is again a Coxeter group with a canonical set of generators which is indexed by
the set of orbits of $\Gamma$ in $I$. Hence, there is a canonical Coxeter diagram $M_{\Gamma}$
over the set of $\Gamma$-orbits in $I$.

The set of stabilized simplices of the  action of $\Gamma$ on $\Sigma_W$ can be canonically identified with
the Coxeter complex associated with $W_{\Gamma}$. In  the sphere representation this corresponds to a fixed lower-dimensional subsphere of the same radius. 



\subsection{Buildings}

A \emph{building of type $(W,S)$} (where $(W,S)$ is a Coxeter system), as introduced by Jacques Tits in the early 60's, is a chamber complex, equipped with a set of simplicial subcomplexes (called \emph{apartments}) isomorphic to the Coxeter complex associated to $(W,S)$, satisfying the three conditions below.
\begin{itemize}
\item[{(B1)}] Every panel is contained in at least three chambers.
\item[{(B2)}] Every two chambers are contained in a common apartment
\item[{(B3)}] For every two apartments there is an isomorphism between these fixing the intersection.
\end{itemize}

If only the last two conditions are fulfilled, then we speak about a \emph{weak building}.

A (weak) building will be \emph{spherical}, \emph{affine}, \emph{irreducible} if its Coxeter system is. It is possible to extend the type function on the simplices of the Coxeter complex to the simplices of the entire building. The \emph{rank} of a building is the rank of the Coxeter system $(W,S)$. 

A \emph{subbuilding} of a spherical building for us is a subset of the geometric realization of it (considering the apartments as spheres), for which the intersections with apartments form convex subsets of a sphere, and such that for each point in this subset there is an apartment containing this point and an antipodal point of the subset. This corresponds to a completely reducible subcomplex in the sense of Serre (see~\cite{Ser:05}).

The simplices containing a given simplex $S$ form again a building, called the \emph{residue} of $S$. The type of this residue can be obtained by deleting the nodes corresponding to $S$ from the Coxeter diagram.

\subsection{Tits diagrams} \label{titsdiagrams}

In \cite{Tit:66} Tits introduced a sort of decorated Coxeter diagrams, which he called the index of a semi-simple algebraic group.
Tits diagrams are combinatorial generalizations of this index. A \emph{spherical Tits diagram} is a triple $(M,\Gamma,A)$
consisting of a spherical Coxeter diagram $M$ (over a set $I$), a subgroup $\Gamma$ of the automorphism group of $M$,
and a $\Gamma$-invariant subset $A$ (called the \emph{anisotropic kernel} of the Tits diagram).
Moreover, it is required 
that each $\Gamma$-orbit $\omega$  in $I \setminus A$
is closed by opposition with respect to the diagram restricted to $A \cup \omega$. In this paper the structure of the group $\Gamma$ is less important. For our purposes it will always suffice to know
the orbits of $\Gamma$ in $I$.
A $\Gamma$-orbit not contained in $A$ is called \emph{isotropic} and the 
\emph{relative rank} of a Tits diagram is the number of its isotropic orbits.

We would like to point out the following observation:
If one deletes some isotropic orbits of a Tits diagram, the remaining diagram is also a Tits diagram.
Thus we can talk of the subdiagrams of relative rank 1 of a Tits diagram.

There is a convenient way to represent a Tits diagram by adding some decoration to the underlying Coxeter diagram: 
The orbits of $\Gamma$ are indicated by putting the nodes
in one orbit close to each other and one encircles the isotropic orbits.

Tits diagrams are useful to describe situations where the fixed point set of a group acting on a building is again a building.
More precisely, given a group $G$ acting on a simplicial building $\Delta$ such that the set of fixed simplices
constitutes a thick building, then the cotype $A$ of a chamber of the fixed building does not depend on the particular chamber
and $(M,\Gamma,A)$ is a Tits diagram (where $M$ is the type of $\Delta$ and $\Gamma$ is the automorphism group of $M$
induced by $G$). In such a situation we say that \emph{the action of the group $G$ admits the Tits diagram $(M,\Gamma,A)$}.



\subsection{The minimal angle of a Tits diagram}\label{section:angle}

In this subsection we will define the \emph{minimal angle of a spherical Tits diagram} $(M,\Gamma,A)$. This will be done in several steps.
If $\Gamma$ is the trivial group and the relative rank is one, then the minimal angle is the angular distance of $i$ with respect to $M$,
where $i$ is the unique index not contained in $A$;
for arbitrary relative rank and trivial $\Gamma$ the minimal angle is defined to be the minimum of the minimal angles in each subdiagram
of relative rank 1.

In order to define the minimal angle in the general case, we consider the diagram $M_{\Gamma}$ and we let $\tilde{A}$
denote the set of $\Gamma$-orbits in $A$. The minimal angle of $(M,\Gamma,A)$
is then defined to be the minimal angle of the Tits-diagram $(M_{\Gamma}, \{ \id \},\tilde{A})$.

Here are three examples:

\begin{center}
\begin{tikzpicture}
\fill (0,0) circle (2pt);
\fill (1,0) circle (2pt);

\draw (0,.15) arc (90:270:.15);
\draw (1,-.15) arc (-90:90:.15);
\draw (0,-.15) -- (1,-.15); 
\draw (0,.15) -- (1,.15); 

\draw (1,.05) -- (0,.05); 
\draw (1,-.05) -- (0,-.05); 

\end{tikzpicture}
\end{center}

In this example $M_{\Gamma}$ is the diagram of type $\mathsf{A}_1$. This is an example of the quasi-split case mentioned in the introduction (so the anisotropic kernel is the empty set) where the minimal
angle is always $\pi$.

\medskip

As a second example, suppose we have the following Tits diagram:

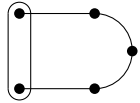
\begin{center}
\begin{tikzpicture}
\fill (0,0) circle (2pt);
\fill (1,1) circle (2pt);
\fill (1,0) circle (2pt);
\fill (0,1) circle (2pt);
\fill (1.5,.5) circle (2pt);

\draw (0,0) -- (1,0); 
\draw (0,1) -- (1,1);

\draw (1,0) arc (-90:90:.5);
\draw (-.15,0) arc (-180:0:.15);
\draw (-.15,1) arc (180:0:.15);
\draw (-.15,0) -- (-.15,1); 
\draw (.15,0) -- (.15,1); 

\end{tikzpicture}\
\end{center}

In this case the diagram $(M_{\Gamma},\{ \id \},\tilde{A})$ is the following:

\begin{center}
\begin{tikzpicture}
\fill (0,0) circle (2pt);
\fill (1,0) circle (2pt);
\fill (2,0) circle (2pt);

\draw (0,0) circle (.15);

\draw (0,0) -- (1,0); 
\draw (1,.05) -- (2,.05); 
\draw (1,-.05) -- (2,-.05); 

\end{tikzpicture}
\end{center}

Here the angular distance at the encircled vertex with respect to $M_{\Gamma}$
is $\pi/2$, which is easily verified. As the Tits diagram is of relative rank 1
it follows that its minimal angle is $\pi/2$.

The third example is a Tits diagram of relative rank 3.
\begin{center}
\begin{tikzpicture}
\foreach \x in {1,...,7}{
\fill (\x,0) circle (2pt);}
\draw (1,0) -- (7,0);
\draw (2,0) circle (.15);
\draw (4,0) circle (.15);
\draw (6,0) circle (.15);
\end{tikzpicture}
\end{center}

Its subdiagrams of relative rank 1 are of the form:

\begin{center}
\begin{tikzpicture}
\foreach \x in {1,...,5}{
\fill (\x,0) circle (2pt);}
\draw (1,0) -- (3,0);
\draw (2,0) circle (.15);
\end{tikzpicture}
\end{center}

Observe that one only needs to consider the connected component containing the isotropic orbit. The minimal angle of this Tits diagram of type $\mathsf{A}_3$ (and so also of the original one of type $\mathsf{A}_7$) is $\pi/2$ (see Section~\ref{section:diag}).

In Section \ref{subsection:angle} further examples of Tits diagrams with minimal angle strictly greater than $\pi/3$
are given.

\subsection{$\R$-buildings}
The notion of $\R$-buildings was introduced by Jacques Tits in~\cite{Tit:86} as a nondiscrete generalization of affine buildings.  Our methods remain valid in this more general case, so we adopt this point of view in our paper.

\subsubsection{Definition}
Let $(\overline{W},S)$ be a spherical Coxeter system of rank $n$ ($n\geq 1$). As mentioned in Section~\ref{section:coxeter}, the group $\overline{W}$ can be realized as a finite reflection group acting on an $n$-dimensional Euclidean space $\mathbb{A}$, which we call the \emph{model space}. Let $W$ be the group of isometries acting on $\A$ generated by $\overline{W}$ and the translations of $\A$.


Consider a pair $(\Lambda,\cF)$ where $\Lambda$ is the set of \emph{points}, and $\cF$ a set of injections (called \emph{charts}) from $\A$ into $\Lambda$. An image of the model space under a chart is called an \emph{apartment}, an image of a root a \emph{half-apartment} and an image of a wall or sector(-face/panel) is called again a \emph{wall} or \emph{sector(-face/panel)}. The pair $(\Lambda,\cF)$ is an \emph{$\R$-building of type $(W,S)$} if the following 6 properties are satisfied (see~\cite[Thm. 1.21]{Par:00}):
\begin{itemize}
 \item[(A1)] If $w\in W$ and $f\in \cF$, then $f \circ w \in \cF$.
 \item[(A2)] If $f,f' \in \cF$, then $X:=f^{-1}(f'(\A))$ is a closed and convex subset of $\A$, and $f|_X = f'\circ w|_X$ for some $w\in W$.
 \item[(A3)] Any two points of $\Lambda$ lie in a common apartment.
\end{itemize}
Due to these last two conditions one can define a distance function $\dd$ on $\Lambda$ such that the injections $\cF$ become isometric embeddings.
\begin{itemize}
 \item[(A4)] Any two sectors $S_1$ and $S_2$ contain subsectors $S_1' \subset S_1$ and $S_2' \subset S_2$ lying in a common apartment.
 \item[(A5)] If three apartments intersect pairwise in half-apartments, then the intersection of all three is nonempty.
\item[(TI)] The distance function $\dd$ is a metric, i.e. it satisfies the triangle inequality.
\end{itemize}


The \emph{dimension} of an $\R$-building is the dimension of its model space. One-dimensional $\R$-buildings are also known as \emph{$\R$-trees}, or shortly trees when no confusion can arise. An $\R$-building is \emph{irreducible} if it does not decompose as direct product of $\R$-buildings.

Different choices for $\cF$ may provide the same metric $\dd$ on $\Lambda$. However there is a unique maximal choice for $\cF$, called the \emph{maximal system of apartments}. One advantage of this choice is that each isometry preserves the set of apartments.

\subsubsection{Global and local structure}\label{section:locglo}
In this section we will  associate (weak) spherical buildings to an $\R$-building $(\Lambda,\cF)$ in two ways.

Two sector-faces are \emph{parallel} if the Hausdorff distance between both is finite. This relation is an equivalence relation due to the triangle inequality. The equivalence classes (named \emph{simplices at infinity}, or the \emph{direction} $F_\infty$ of a sector-face $F$) form a spherical building $\Lambda_\infty$ of type $(\overline{W},S)$ called the \emph{building at infinity} of the $\R$-building $(\Lambda,\cF)$. The chambers of this building are the equivalence classes of parallel sectors. There is a bijective correspondence between the apartments of $(\Lambda,\cF)$ and the apartments of $\Lambda_\infty$. The building at infinity corresponding to the maximal system of apartments is called the \emph{complete building at infinity}. Another way to define the complete building at infinity is using the boundary at infinity of the metric space (see~\cite[4.2.1]{Kle-Lee:97}). The rank of the building at infinity equals the dimension of the $\R$-building. 

\begin{remark}
\emph{We will always assume that the building at infinity is actually a building, rather than a weak building. One can always reduce the second case to the first case using~\cite[4.9]{Kle-Lee:97}.}
\end{remark}

Two sector-faces are \emph{asymptotic} if they have a sector-face of the same dimension as the original two in common. This is also an equivalence relation. Asymptotic sector-faces are necessarily parallel, the converse is only true for sectors (see~\cite[Cor.~1.6]{Par:00}).

One can also define local equivalences. The apex of a sector(-face) is called its \emph{base point}, or we say that the sector(-face) is \emph{based} at this point. Let $\alpha$ be a point of $\Lambda$, and $F,F'$ two sector-faces based at $\alpha$. Then these two sector-faces will \emph{locally coincide} if their intersection is a neighborhood of $\alpha$ in both $F$ and $F'$. This relation forms an equivalence relation defining \emph{germs of sector-faces} as equivalence classes. These germs form a (weak) spherical building  $\Lambda_\alpha$ of type $(\overline{W},S)$, called (again) the \emph{residue} at $\alpha$.  

A detailed study of $\R$-buildings can be found in the article~\cite{Par:00}.
We collect several results of that paper in order to have them later  available for reference:

\begin{lemma}[\cite{Par:00}, Prop.~1.16]\label{p16}
Let $S$ and $S'$ be two sectors, then there exists an apartment containing two sectors with the same germ as $S$ and $S'$ respectively.
\end{lemma}

 \begin{lemma}[\cite{Par:00}, Prop.~1.8]\label{p8}
Let $C$ be a chamber of the building at infinity $\Lambda_\infty$ and $S$ a sector based at $\alpha \in
\Lambda$. Then there exists an apartment $\Sigma$ containing an element of the germ of $S$ and such that $\Sigma_\infty$ contains $C$. 
 \end{lemma}

 \begin{cor}[\cite{Par:00}, Cor.~1.9]
Let $\alpha$ be a point of $\Lambda$ and $F$ a simplex of the building at infinity. Then there is a unique
sector-face based at $\alpha$ with direction $F$.
 \end{cor}

The unique sector-face of the previous corollary will be denoted by $F_\alpha$. We will often use the subscript to indicate the base point of the sector-face, or use $\infty$ as subscript to denote the direction of the sector-face. A consequence of this corollary is the existence of canonical, type-preserving epimorphisms from the building at infinity to the residues.


\begin{lemma}[\cite{Par:00}, Prop.~1.12]\label{p12}
If two sectors $S_\alpha$, $S'_\alpha$ have germs forming opposite chambers of $\Lambda_\alpha$, then there exists a unique apartment containing both.
\end{lemma}

\subsubsection{Classification results and the little projective group} \label{section:classification}
Irreducible spherical buildings of rank at least three satisfy a powerful transitivity property, called the \emph{Moufang property}.  Moufang buildings of rank at least two are classified and are shown to arise from certain classical, algebraic and mixed groups over (skew) fields. For details see~\cite{Tit:74},~\cite{Tit-Wei:02} and~\cite{Wei:03}. 

An important subgroup of the automorphism group of such a building, and in particular the group the building arose from, is the group generated by all the root groups of the building. This group is denoted by $G^\dagger$ in~\cite[Lem.~11.20]{Wei:03}. In the rank one and two cases this group is often called the \emph{little projective group} (see for example~\cite{Med-Seg-Ten:08} and~\cite[Def. 4.4.4]{Mal:98}); we will use this last denomination regardless of rank.

One can classify irreducible $\R$-buildings of dimension at least three, see~\cite{Tit:86} and also~\cite{Wei:09}, using the fact that their buildings at infinity are Moufang. These $\R$-buildings correspond to Moufang buildings with the additional information of a valuation on the underlying (skew) field satisfying certain properties. The buildings that arise in this way are called \emph{Bruhat-Tits buildings}. We will be mainly interested in the algebraic group case, and the mixed group case to a lesser extent. 

%
%
%
%

\subsubsection{Bruhat-Tits buildings for algebraic groups and the existence problem}\label{section:bruhattits}

Let $G$ be a connected reductive algebraic group defined over a field $K$ with a valuation $\nu$. To these data one wants to  associate a Bruhat-Tits building $I_{K,\nu}(G)$  on which the group $G(K)$ acts (see~\cite[\S 2]{rousseauthesis} for details).

It is conjectured in~\cite[p. 173]{Tit:86} that a sufficient condition for the existence of the Bruhat-Tits building $I_{K,\nu}(G)$ is that the valuation is complete. A positive answer is known in various cases, for example the case where $G$ is a classical group (see \cite{Bru-Tit:72}) or when $G$ is quasi-split over $K$ (see~\cite[\S 4]{Bru-Tit:84}).  The main goal of this two-part paper is to provide a positive answer for all algebraic cases.

One basic strategy to attack the conjecture is to consider a finite Galois extension $L$ of $K$ such that the algebraic group $G(L)$ is quasi-split (which exists by~\cite[5.1.4]{Bru-Tit:84}). The valuation $\nu$ extends uniquely to a complete valuation $\omega$ of $L$ by~\cite[Thm. 7.1.1]{Cas:86}.  As mentioned above, one can now consider the Bruhat-Tits building $I_{L,\omega}(G)$. The Galois group $\mbox{Gal}(L/K)$ acts  on $I_{L,\omega}(G)$, again by~\cite[5.1.4]{Bru-Tit:84}. The fixed structure of $\mbox{Gal}(L/K)$ on the complete building at infinity is the spherical building $\Delta_K(G)$, moreover this action admits a Tits diagram for which the possibilities are listed in~\cite{Tit:66}.


\subsubsection{Bounded isometries and fixed point theorems}\label{section:bounded}
Let $(X,\dd)$ be a \emph{geodesic metric space}, which is by definition a metric space where every two points can be joined by a \emph{geodesic segment}, which is an isometrically embedded line segment. Let $\Delta$ be a triangle in $X$ with geodesic segments as sides. As $X$ is a metric space, one can construct in the Euclidean plane a triangle $\Delta'$ with the same side lengths as $\Delta$. If the distances between points of the sides of $\Delta$ are less than or equal to the distances between the corresponding points on $\Delta'$, we say that $(X,\dd)$ is a CAT(0)-\emph{space}.

The metric spaces formed by $\R$-buildings, as well as their completions and convex subsets, are CAT(0)-spaces (see \cite[Lem.~3.2.1]{Bru-Tit:72}). A  first fixed point theorem for CAT(0)-spaces we will need is the following:

\begin{theorem}[\cite{Bru-Tit:72}, Prop.~3.2.4]\label{bt}
A nonempty bounded subset of a complete CAT(0)-space $X$ has an unique `center'.
\end{theorem}

An isometry or a group of isometries of a CAT(0)-space is called \emph{bounded} if it admits a bounded orbit. Note that if some finite index subgroup of a group $G$ of isometries is bounded, then so is $G$ itself. The following direct corollary of the above lemma for such groups is known as the Bruhat-Tits fixed point theorem.

\begin{cor}\label{bt2}
If $G$ is a bounded group of isometries of a complete CAT(0)-space $(X,\dd)$, then $G$ fixes some point in $X$.
\end{cor}

We mention that affine buildings are always metrically complete. On the other hand nondiscrete $\R$-buildings are often noncomplete (see for instance~\cite{Mar-etal:*}). This does not pose a major problem for our purposes however, by the following remark.

\begin{remark}\label{rmk} \rm
Throughout this paper we will apply Theorem~\ref{bt} and its corollary a couple of times. Note that, when the $\R$-building is not metrically complete, we cannot apply these results directly. However, the completion of the metric space  on $\Lambda$ is a complete CAT(0)-space, so there exists a center in the completion. One can work with points in the completion as if they were points of an $\R$-building; this is because one can always embed the $\R$-building $(\Lambda,\cF)$ in a metrically complete $\R$-building  of the same type (see for example~\cite[Lem.~4.4]{Str:*}). For these reasons, when we use the notions `center' or `fixed point', we implicitly allow them to be in the completion.
\end{remark}

The next type of fixed point theorems deals with groups consisting solely of bounded isometries.

\begin{theorem}[\cite{Mor-Sha:85}, Prop.~II.2.15]\label{thm:trees}
If a finitely generated group $G$ of isometries of an $\R$-tree consists solely of bounded elements, then $G$ itself is bounded.
\end{theorem}

\begin{theorem}[\cite{Par:03}, Cor.~3]\label{thm:par}
Let $G$ be a connected reductive algebraic group defined and quasi-split over a field $K$ with complete valuation $\nu$. If a finitely generated subgroup $H$ of $G(K)$ consists solely of bounded elements w.r.t.~their action on $I_{K,\nu}(G)$, then $H$ is bounded.
\end{theorem} 

Provided that one has a positive (partial) answer to the existence problem mentioned in Section~\ref{section:bruhattits}, one can use the following theorem to generalize Theorem~\ref{thm:par}.


\begin{theorem}\label{theorem:embed}
Let $G$ be a connected reductive algebraic group defined over a field $K$ with valuation $\nu$, and let $(L,\omega)$ be an extension of $(K,\nu)$. If the Bruhat-Tits buildings $I_{K,\nu}(G)$ and $I_{L',\omega'}(G)$ exist for every extension $(L',\omega')$ of $(K,\nu)$, then $I_{K,\nu}(G)$ embeds $G(K)$-equivariantly in the completion of $I_{L,\omega}(G)$.
\end{theorem}
\proof
This essentially follows from~\cite[Thm. 5.3.3]{rousseauthesis}. This theorem, however, needs a second assumption concerning the metrical completeness of a  Bruhat-Tits building $I_{L,\omega}(Z')$, where $Z$ is the anisotropic kernel. This Bruhat-Tits building is isometric to the $\R$-building $T(B)$ obtained from Section~\ref{section:gentrees}, where $B$ is a certain stabilized subspace. This extra condition comes into play in part 2.c of the proof~\cite[5.3.5]{rousseauthesis} of the aforementioned theorem. It is used to find a fixed point of a bounded subgroup acting on $I_{L,\omega}(Z')$. Removing the completeness condition yields a point in the completion instead, which corresponds to a subspace with the same direction as $B$ in the completion of $I_{L,\omega}(G)$ (note that by Remark~\ref{rmk} we can treat points, as well as subspaces, in the completion of an $\R$-building in the same way as those in the building itself). Following the proof further down gives the same conclusion as the original theorem, with the added possibility that $I_{K,\nu}$ is embedded in the completion of $I_{L,\omega}$ instead of $I_{L,\omega}$ itself.
\qed 

If one is only interested in the existence problem for exceptional algebraic groups, then Theorem~\ref{theorem:embed} is not needed, see Remark~\ref{rem:steps}.

\begin{cor}\label{cor:embedcor}
Let $G$ be a connected reductive algebraic group defined over a field $K$ with complete valuation $\nu$. Suppose that the Bruhat-Tits buildings $I_{K,\nu}(G)$ and $I_{L,\omega}(G)$ exist for every extension $(L,\omega)$ of $(K,\nu)$. If a finitely generated subgroup $H$ of $G(K)$ consists solely of bounded elements w.r.t.~their action on $I_{K,\nu}(G)$, then $H$ is bounded.
\end{cor}
\proof
Under these assumptions Theorem~\ref{theorem:embed} allows us to  embed the Bruhat-Tits buildings $I_{K,\nu}(G)$ $G(K)$-equivariantly in the completion of $I_{L,\omega}(G)$, where $L$ is  a finite Galois extension such that the algebraic group $G(L)$ is quasi-split (see Section~\ref{section:bruhattits}). This allows us to interpret $H$ as consisting solely of bounded elements w.r.t.~their action on the completion of $I_{L,\omega}(G)$. By Theorem~\ref{thm:par} $H$ is bounded w.r.t.~its action on the completion of $I_{L,\omega}(G)$ and so also w.r.t.~its action on $I_{K,\nu}(G)$. \qed

The following remark will be of interest in the proof of our Main Result in Section~\ref{section:proof} when we apply Corollary~\ref{cor:embedcor}.

\begin{remark}\label{rem:cases} \rm
If one considers the possible residues of non-trivial simplices of spherical buildings of rank at least three which are not associated to classical non-algebraic groups, then by considering~\cite[\S 41]{Tit-Wei:02} and~\cite[\S 12]{Wei:03} one of the following three possibilities occurs.
\begin{enumerate}
\item
The residue is a spherical building associated to a group which is both classical and algebraic. The existence problem for the corresponding Bruhat-Tits buildings is solved in~\cite{Bru-Tit:72}. Hence Corollary~\ref{cor:embedcor} applies.
\item
The residue is the spherical building (associated to a group of so-called mixed type) $\mathsf{C}_l(K,F,1)$, where $K$ is a field of characteristic 2, $F$ is a subfield of $K$ such that $K^2 \subset F \subset K$ and $l$ equals 2 or 3. These buildings, and their associated Bruhat-Tits buildings embed $G^\dagger$ equivariantly (with $G^\dagger$ being the little projective group) respectively in the buildings $\mathsf{B}_l(K,K,q)$ and their associated Bruhat-Tits buildings, where $q$ is the quadratic form from $K$ to itself by taking squares. The  latter buildings  correspond with split algebraic groups. By this equivariant embedding one proves an analogue of Theorem~\ref{thm:par} similarly as Corollary~\ref{cor:embedcor}.
\item
The residue is an exceptional algebraic building associated to the following Tits diagram. 
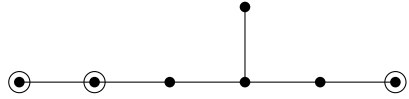
\begin{center}
\begin{tikzpicture}
\foreach \x in {1,...,6}{
\fill (\x,0) circle (2pt);}
\draw (2,0) circle (4pt);
\draw (1,0) circle (4pt);
\draw (6,0) circle (4pt);

\fill (4,1) circle (2pt);

\draw (1,0) -- (6,0); 
\draw (4,0) -- (4,1); 

\end{tikzpicture}
\end{center}
The minimal angle of Tits diagram is $\pi/2$ (see Section~\ref{subsection:angle}). Note that if the Main Result (see Section~\ref{mainresult}) holds for Bruhat-Tits buildings associated to a quasi-split algebraic group, then this Main Result implies that there is again an equivariant embedding so we can again prove an analogue of Theorem~\ref{thm:par}.
\end{enumerate}
\end{remark}

%
%
%

The next lemma allows us to extend the above results to the infinitely generated groups, where an extra possibility occurs.

\begin{lemma}\label{lem:cap}
Let $G$ be a group of isometries of an Euclidean building $(\Lambda,\cF)$. If every finitely generated subgroup of $G$ is bounded, then $G$ is bounded or there is a non-empty simplex of the complete building at infinity stabilized by $G$ and its centralizer in the isometry group of $(\Lambda,\cF)$.
\end{lemma}
\proof
The fixed sets of the finitely generated subgroups of $G$ form a filtering family of closed convex subsets of the metric completion $\overline{\Lambda}$. Such a filtering family has a common element in $\overline{\Lambda}$ or it fixes a subset of intrinsic radius at most $\pi/2$ in its boundary by~\cite[Thm.~1.1]{Cap-Lyt:09}. The first possibility implies that $G$ is bounded, the second that there is a non-empty simplex of the complete building at infinity stabilized by $G$ and its centralizer in the isometry group of $(\Lambda,\cF)$ by~\cite[Prop.~1.4]{Bal-Lyt:05}.\qed

The next theorem describes the nature of unbounded isometries. The \emph{translation length} of an isometry $g$ is the infimum $\inf\{\dd(x,g(x))\}$.

\begin{theorem}[\cite{Par:00}, \S4] \label{thm:para}
Let $g$ be an isometry of an $\R$-building $(\Delta,\cF)$. Then $g$ fixes some point or translates a geodesic line in the completion $\overline{\Delta}$. The union of all such geodesic lines is exactly the set of points with displacement the translation length.
\end{theorem}

Isometries which translate a geodesic line in a non-trivial way are called \emph{hyperbolic} isometries.

\section{Statement of the main result} \label{mainresult}

We are now able to state the main result.

\textbf{Main Result.} \emph{
Let $(\Lambda,\cF)$ be an irreducible $\R$-building with a maximal system of apartments, and let $G$ be a bounded group of isometries acting on $\Lambda$. Suppose that the following holds:
\begin{itemize}
\item The action of $G$ on $\Lambda_\infty$ admits a Tits diagram whose minimal angle is strictly greater than $\pi/3$.
\item $(\Lambda,\cF)$ is not a Bruhat-Tits building arising from a classical non-algebraic group.
\end{itemize}
Then the fixed structure of $G$ in the completion of $\Lambda$ contains an $\R$-building $(\Lambda',\cF')$ such that
$\Lambda'_\infty$ is precisely the fixed point set of $G$ in $\Lambda_\infty$, and that the isometries of $\Lambda$ centralizing $G$ (and thus acting on the spherical building $\Lambda'_\infty$) also act on $(\Lambda',\cF')$ in a canonical way.
}


An application of the main result will be discussed in Section~\ref{exis}. 

\begin{remark}\rm
The restriction involving classical non-algebraic group stems from the algebraic nature of Theorem~\ref{thm:par}. In particular, if one can extend this Theorem to any $\R$-building, then this restriction can be lifted from the main result.
\end{remark}

\section{Generalizations of wall- and panel-trees}\label{section:gentrees}
This section is devoted on how to derive new (lower-dimensional) $\R$-buildings from a given $\R$-building. These constructions generalize the so-called wall- and panel-trees (see for instance~\cite[Prop.~4]{Tit:86}).

These more general constructions are not new (see for example~\cite[Prop.~4.8.1]{Kle-Lee:97} or~\cite[Prop.~8.1.5]{Cha:09}), however we will include the constructions and the proofs of the basic facts concerning them in order to ensure generality and for future reference.

\subsection{Construction}
Let $(\Lambda,\cF)$ be an $\R$-building with apartment set $\cA$. Let $S_\infty$ and $S'_\infty$ be two opposite non-maximal simplices at infinity; the residue at infinity  of $S_\infty$ is a spherical building $(\Lambda_\infty)_{S_\infty}$. 

Let $B$ be the smallest convex subcomplex of the building at infinity containing the simplices $S_\infty$ and $S'_\infty$. This complex forms a subsphere of the same dimension as $S_\infty$. An affine subspace $M$ of an apartment of the $\R$-building with this subcomplex $B$ at infinity will be referred to as a \emph{subspace with direction $B$} (denoted by $M_\infty =B$). By this we mean that the direction of each sector-face contained in this subspace is contained in $B$, and conversely that each sector-face with direction in $B$ based at a point of the subspace is completely contained in the subspace.

Let $T(B)$ be the set of all subspaces $M$ of the $\R$-building with $M_\infty = B$, and $T(S_\infty)$ the set of all asymptotic classes of sector-faces in the parallel class $S_\infty$. 


We now define charts on $T(B)$ and $T(S_\infty)$: Choose a subspace $M$ and a sector-face $D$ of the model space, such that there exists some chart $f \in \cF$ with $f(M)_\infty=B$ and $f(D) \in S_\infty$. One can identify the model space $\A$ with the product $\R^m \times M$ for a certain $m \in \mathbb{N}$. For every chart $g \in \cF$ of the $\R$-building $(\Lambda,\cF)$ with $g(M)_\infty =B$ or  $g(D) \in S_\infty$, one defines a chart $g'$ on $T(B)$ or $T(S_\infty)$  respectively as follows: if $r \in \R^m$, then $g'(r)$ is the subspace $g(\{r\} \times M)$, or the asymptotic class containing $g(\{r\} \times D)$ respectively.





\subsection{Proof}
Before proving that these two constructions yield $\R$-buildings we show that they are equivalent. For this we need a few lemmas. When we mention `subsector-face', we mean sector-faces which are subsets of the other sector-face and of the same dimension.
\begin{lemma}\label{lemma:inapp}
Let $S_\alpha$ and $S_\beta$ be two sector-faces based at respectively $\alpha,\beta \in \Lambda$ with the same direction $S_\infty$. Then there exists an apartment containing subsector-faces of both.
\end{lemma}
\proof
Because the two sector-faces are parallel it is sufficient to prove that  the point $\beta$ and some subsector-face of $S_\alpha$ lie in a common apartment. Embed $S_\alpha$ in some apartment $\Sigma$. Let $T_\infty$ be the simplex of $\Sigma_\infty$ opposite to $S_\infty$. Given a point $\gamma$ of $S_\alpha$ let $K_\gamma$ be the union of all sectors in $\Sigma$ containing $T_\gamma$. Let $B$ be an open ball in $\Sigma$ with center $\alpha$ and radius strictly greater than $\dd(\alpha,\beta)$. We fix $\gamma$ such that $K_\gamma$ contains $B$. 

Let $\Sigma'$ be an apartment containing $\beta$ and a germ of $S_\gamma$ (possible by Lemma~\ref{p16}). Let $R_\gamma$ be a sector in $\Sigma$ based at $\gamma$ containing $S_\gamma$, and $R'_\gamma$ a sector in $\Sigma'$ based at $\gamma$ containing $\beta$. The canonical retraction $r$ of $\Lambda$ on $\Sigma$ centered at $R_\gamma$ (see~\cite[Prop.~1.17]{Par:00}) maps the sector $R'_\gamma$ to a sector $C_\gamma$ in $\Sigma$. As this retraction is distance nonincreasing, the point $r(\beta)$ lies in the ball $B$, so the sector $C_\gamma$ is contained in $K_\gamma$. In particular this implies that the germ of $C_\gamma$ contains a simplex opposite to the germ of $S_\gamma$. By the way the retraction is defined it follows that the germ of $R'_\gamma$ contains a simplex with the same property.

From~\cite[Prop.~9.9]{Wei:03} it follows that there exists a chamber in the residue at $\gamma$ containing the germ of the sector-face $S_\gamma$ and opposite the germ of $R'_\gamma$. Using Lemma~\ref{p8} one can find a sector $D_\gamma$ with this chamber as germ and containing $S_\gamma$. We conclude by Lemma~\ref{p12} that there exists an apartment containing $\beta$ and a subsector-face of $S_\alpha$. \qed



\begin{lemma}\label{lemma:helpie}
Let $S_\alpha$ be a sector-face based at $\alpha \in \Lambda$ and $S'_\infty$ a simplex at infinity opposite to $S_\infty$. Then there exists a unique subspace (of the same dimension as $S_\alpha$) containing both $S'_\infty$ at infinity and a subsector-face of $S_\alpha$.
\end{lemma}
\proof
Let $B$ be some subspace containing both $S'_\infty$ and $S_\infty$ at infinity. Let $\beta$ be a point of this subspace. By the previous lemma there exists an apartment $\Sigma$ containing subsector-faces of both $S_\alpha$ and $S_\beta$. In particular there exists a sector $C_\gamma$ based at some point $\gamma$ of $S_\beta$ containing a subsector-face of $S_\alpha$. The germ of this sector is opposite to some germ of a sector $D_\gamma$ containing $S'_\gamma$. It is clear that the unique apartment provided by Lemma~\ref{p12} containing $C_\gamma$ and $D_\gamma$ contains a desired subspace. 

Uniqueness follows directly as distinct subspaces with the same direction are disjoint (by the definition of a subspace).
\qed

The above lemma makes clear that the sets of points of the two constructions are in a one-to-one correspondence with each other. An apartment of $T(B)$ is easily seen to give rise to an apartment of $T(S_\infty)$. Conversely, for $T(S_\infty)$, one sees that all apartments containing $S_\infty$ and an apartment in the residue of $S_\infty$ at infinity (in the sense of spherical buildings), correspond to one apartment of $T(B)$, establishing a one-to-one correspondence for the apartments in both constructions.

We now verify (A1)-(A5) and (TI) for both constructions. The above discussion implies that we can choose which construction to verify the condition for. 

\begin{itemize}
\item[{(A1),(A2)}]
Directly from the corresponding conditions of the original building and the construction of $T(B)$.
\item[{(A3)}] 
From Lemma~\ref{lemma:inapp} using the construction of $T(S_\infty)$.
\item[{(A4)}]
Notice that sectors in the second construction are in fact sectors of the original building, so (A4) for the original building gives us subsectors contained in one apartment, which yields an apartment in the construction of $T(S_\infty)$.
\item[{(A5)}]
From the construction of $T(B)$ and (A5) for the original building.
\item[{(TI)}]
The Hausdorff metric (which satisfies the triangle inequality) yields for $T(B)$ the natural distance function on the apartments of $T(B)$. Hence (TI) is satisfied.
\end{itemize}

As the conditions are all verified, we have proved that these constructions yield $\R$-buildings. Finally we remark that the building at infinity of $T(S_\infty)$ will be the residue of $S_\infty$ in the building at infinity (which is seen by considering sectors in $T(S_\infty)$).



\section{Proof of the main result}\label{section:proof}

Suppose we are given a situation as described in the statement of the main result. As the action of the group $G$ at infinity admits some Tits diagram $(M,\Gamma,A)$, its fixed structure forms a spherical building of type $(\overline{W}',S')$.

We prove the main result using an induction step. In Step one we go through the proof assuming that $(\Lambda,\cF) $ is a Bruhat-Tits building associated to a quasi-split group. After that in Step two we repeat the proof dealing with general Tits diagrams.


\begin{remark}\label{rem:steps} \rm
If one is only interested in the existence problem for exceptional algebraic groups, then one does not need to use step two, and in particular Theorem~\ref{theorem:embed} which is used therein.
 \end{remark}

\subsection{Stabilized subspaces}
Let $S_\infty$ and $S'_\infty$ be maximal stabilized and opposite simplices at infinity (note that their types are the same), and let $B$ be the smallest convex subcomplex of $\Lambda'_\infty$ containing both.

\begin{lemma}\label{exist1}
There exists at least one subspace with direction $B$ and stabilized by $G$.
\end{lemma}
\proof
As $G$ is bounded, it fixes at least one point $\alpha$ by Corollary~\ref{bt2}. The sector-face $S_\alpha$ based at this point with direction $S_\infty$ is stabilized by $G$, so the `center geodesic ray' of this sector-face is fixed by it. A consequence is that each subsector-face of $S_\alpha$ contains a fixed point.

It follows that the subspace containing a subsector-face of $S_\alpha$ and the simplex $S'_\infty$ at infinity constructed by Lemma~\ref{lemma:helpie} contains a fixed point. As $S_\infty$ and $S'_\infty$ are both stabilized, we conclude that this subspace is stabilized as well.
\qed

Note that the dimension of such a subspace equals the number of encircled nodes in the Tits diagram $(M,\Gamma,A)$. An isometric embedding of a Euclidean space in an apartment of the $\R$-building is called a \emph{flat}.

\begin{lemma}\label{fixp}
The fixed points of $G$ in a stabilized subspace with direction $B$ form a flat with dimension the relative rank of the Tits diagram.
\end{lemma}
\proof
Let $M$ be such a stabilized subspace. The group $G$ stabilizes this subspace, hence $G$ acts as a bounded group of isometries of $M$. As $M$ is evidently a complete CAT(0)-space, $G$ fixes at least one point $\alpha$ of $M$ by Corollary~\ref{bt2}. 

Let $X$ be the set of points in $M$ fixed by $G$. If a point of $M$ is spanned by points of $X$, then it is fixed as well, hence $X$ is a flat. The assertion about the dimension of the flat follows from the action of $G$ on the stabilized sector-face $S_\alpha$.
\qed


By Lemma~\ref{exist1} we know that there is at least one stabilized subspace with $B$ at infinity; consider all such stabilized subspaces. Our goal is now to fix one such subspace for each possible choice of $S_\infty$ and $S'_\infty$. If these simplices are chambers, then there is only one such subspace (the unique apartment containing both), so suppose that they are not chambers. These subspaces then form a set $F$ of points of the $\R$-building $T(B)$. The set $F$ is closed, convex, nonempty and contains no points in its (visual) boundary (because of the maximality of $S_\infty$ and $S'_\infty$). In particular it does not contain geodesic lines.

\begin{lemma}\label{lemma:nohyp}
There is no hyperbolic isometry of $T(B)$ stabilizing $F$.
\end{lemma}
\proof
We work in the completion $\overline{T(B)}$. Let $g$ be such a hyperbolic isometry. Then $g$ translates some geodesic line $L$ by Theorem~\ref{thm:para}. Remark that such a line is disjoint from $\overline{F}$. Let $C$ be a closed line segment of $L$ such that its orbit under $g$ covers $L$. Each point of $C$ has a (unique) closest point in $\overline{F}$ (see~\cite[Prop.~II.2.4]{Bri-Hae:00}). As $C$ is compact, one hence can find an $\alpha \in \overline{F}$ and $\beta \in C$ such that $\dd(\alpha,\beta)$ is minimal amongst all choices of $\alpha$ and $\beta$. The line segment $L'$ between $\alpha$ and $\alpha^g$ is contained in $\overline{F}$ by convexity. Note that $\dd(\alpha,\beta)$ is the minimal distance between points of $L$ and $L'$. This implies that all of the angles of the quadrangle with corners $\alpha$, $\alpha^g$, $\beta^g$ and $\beta$ (in order) are at least $\pi/2$. By~\cite[Prop.~II.2.11]{Bri-Hae:00} one concludes that the convex hull of these corners is isometric to a rectangle in the Euclidean plane, and hence that $\dd(\alpha,\alpha^g) = \dd(\beta,\beta^g)$. This contradicts Theorem~\ref{thm:para}.
\qed

By the results of Section~\ref{section:gentrees} we have isometries from $T(S_\infty)$ to $T(B)$ and from $T(B)$ to $T(S'_\infty)$. By combining these isometries along paths of subsequent opposite maximal stabilized simplices starting and ending at $S_\infty$ we obtain a group $H$ of isometries of $T(S_\infty)$. This group is centralized by the induced action of $G$ on $T(S_\infty)$. Note that if we consider only paths of even length we obtain a normal subgroup $H' \lhd H$ of index at most 2. 

We now claim that $H$ fixes a point. Lemma~\ref{lemma:nohyp} asserts that $H$ only contains bounded isometries.

 If $T(S_\infty)$ is a tree we are done by Theorem~\ref{thm:trees} and Lemma~\ref{lem:cap} (the existence of a non-empty simplex at infinity fixed by the centralizer of $H$ would imply that $G$ fixes it as well, contradicting the maximality of $S_\infty$). 

If $T(S_\infty)$ is not a tree and we are in Step one, so $(\Lambda,\cF)$ is associated to a quasi-split group, then $T(S_\infty)$ is also associated to a quasi-split group.  Note also that then $(\Lambda,\cF)$ is an $\R$-building of dimension at least three. Repeating the discussion mutatis mutandis of Section 2 of~\cite{Kna:88} one shows that the group $H'$ is a subgroup of the stabilizer of $S_\infty$ in the little projective group of the Moufang building at infinity. This allows us to apply Section~\ref{section:classification} and Theorem~\ref{thm:par}, concluding that $H'$ is bounded and fixes some point. As $H'$ is of index at most 2 in $H$, this holds for $H$ as well.

Now suppose that $T(S_\infty)$ is not a tree and we are in Step two. Note that then $(\Lambda,\cF)$ is an $\R$-building of dimension at least three. From the case study made in Remark~\ref{rem:cases} and the results obtained in Step one, it follows that in each case Theorem~\ref{thm:par} or an analogue of it holds. Hence we can proceed as in Step one and conclude that $H$ is bounded and fixes a point.


By projecting this fixed point onto $F$ (again using~\cite[Prop.~II.2.4]{Bri-Hae:00}), we obtain a fixed point of $H$ in $F$. By using the isometries obtained in Section~\ref{section:gentrees} this corresponds to the choice of a unique point in the $\R$-buildings $T(S_\infty)$, $T(S'_\infty)$ and $T(B)$ which is stable under the canonical isometries between these, for every possible $S_\infty$, $S'_\infty$ and $B$. With this unique point there corresponds a stabilized subspace with direction $B$ (again for every possible $B$). So we arrived at our goal of fixing a stablized subspace with direction $B$.

We call this subspace the \emph{middle stabilized subspace with direction $B$}.  The fixed part of this subspace (as described in Lemma~\ref{fixp}) is called the~\emph{middle fixed flat with direction $B$}. 

The structure at infinity of such a middle fixed flat corresponds to an apartment of the spherical building formed by the fixed point structure at infinity. In particular the boundary (which is a sphere) of this flat  allows for the natural action of the Coxeter group $\overline{W}'$ on it. Combining this with the translations acting on the flat, one can define charts from the model space $\A'$ defined from $(\overline{W}',S')$ to these middle fixed flats.

\begin{remark} \rm \label{rem:applibound}
Note that if $F$ is bounded then one can use the center of this set (see Theorem~\ref{bt}) as the choice of a fixed point. This explains the use of `middle'. For the application concerning Galois descent for algebraic groups this is always the case, see~\cite[Prop.~5.2.1]{rousseauthesis}.
\end{remark}


\subsection{Constructing $(\Lambda',\cF')$}\label{section:subbie}
Let $\Lambda'$ be the union of the points of all middle fixed flats and $\cF'$ be the set of all charts to these fixed middle flats as discussed in the previous section. So we can talk about apartments of $\Lambda'$. One easily sees that  Conditions (A1) and (TI) both are satisfied. The closed and convex part of Condition (A2) is also directly satisfied.


From the way we chose the middle stabilized subspace, and in particular its interplay with the canonical isometries between the various derived $\R$-buildings, it follows that if two apartments of $\Lambda'$ share a maximal stabilized simplex at infinity, then the corresponding sector-faces in both apartments are asymptotic, or by definition, they share a common subsector-face. Condition (A4) now follows from (B2) which states that every two chambers of the fixed building at infinity lie in a common apartment.

The same reasoning combined with the convexity from (A2) shows that if two apartments of $\Lambda'$ share a half-apartment at infinity, then the apartments themselves share a half-apartment. In order to be able to prove Conditions (A3), (A5) and the last part of Condition (A2) we will need the extra assumption that the minimal angle of the Tits diagram is strictly greater than $\pi/3$. For a discussion on when this occurs see Section~\ref{subsection:angle}.


\subsubsection{Condition (A5)}

Assume that there exist three apartments of $\Lambda'$ pairwise sharing a half-apartment, while the intersection of all three is nonempty. Such a configuration we call a \emph{triangle configuration}. If the relative rank of the Tits diagram equals the rank of the Coxeter system $(W,S)$ then this is  impossible due to Condition (A5) for the original building, so suppose this is not the case. 

The three corresponding apartments at infinity share some simplex (which is the trivial simplex if the Tits diagram is of relative rank one). Using the construction from Section~\ref{section:gentrees} on such a shared simplex $S$ of maximal dimension, one obtains a triangle configuration in $T(S)$ of flats isometric to the real line (so these flats are geodesic lines). This reduction to $T(S)$ corresponds with removing all isotropic orbits from the Tits diagram $(M,\Gamma,A)$ but one, resulting in a subdiagram $(M',\Gamma\vert_{M'},A)$ of relative rank one. Let $(W',S')$ be the Coxeter system associated to $M'$.

Let $\alpha$ be a corner of the triangle formed. As $\alpha$ and $S$ are both stabilized by $G$, one can consider the group action of $G$ on the residue $T(S)_\alpha$. This residue is a (weak) building of type $(W',S')$ and the diagram automorphisms induced here are the same as at the building at infinity. This happens because there is a type-preserving epimorphism from the building at infinity to this residue (see Section~\ref{section:locglo}).

The two sides of the triangle meeting in $\alpha$ correspond to stabilized simplices of the group action. Following Sections~\ref{section:diag} and~\ref{section:angle}, these stabilized simplices correspond to 2 vertices of a certain type of the Coxeter diagram $M'_\Gamma$, where the type corresponds to the isotropic orbit of the Tits diagram $(M',\Gamma\vert_{M'},A)$. 

By the assumption on the minimal angle, one has that the angle between these two vertices, and so also of the corner $\alpha$ of the triangle in $T(S)$, is strictly greater than $\pi/3$. But this contradicts the fact that the sum of the angles of a triangle in a CAT(0)-space is less than or equal to $\pi$ (see~\cite[Prop.~II.1.7(4)]{Bri-Hae:00}). This concludes the proof of Condition (A5).



\subsubsection{Conditions (A2) and (A3)}

The only conditions one still has to verify are part of Condition (A2) and Condition (A3). The strategy we apply here is the same as used by Petra N. Schwer and the second author in~\cite{Sch-Str:*}. This section only uses Condition (A5) and no hypothesis on the minimal angle.

One can define sectors, sector-faces and germs of sector-faces in $(\Lambda',\cF')$ as usual. Hence we may define residues; namely simplicial complexes $\Lambda'_\alpha$ for points $\alpha \in \Lambda'$. We now show that these are again (weak) buildings: The proof of Lemma~\ref{p8} remains valid in this setting (see~\cite{Par:00} for the proof), so for every two simplices in a $\Lambda'_\alpha$ one can find an apartment of $\Lambda'$ such that these two simplices appear as germs. This proves Condition (B2). Condition (B3) for $\Lambda'_\alpha$ follows from the same condition for the (weak) building $\Lambda_\alpha$. Hence $\Lambda'_\alpha$ is indeed a (weak) spherical building. We mention that the proof of Lemma~\ref{p12} remains valid as well.

%
%
%

We now deal with the second part of (A2). As we are dealing with isomorphically embedded Euclidean spaces in a metric space, it follows that one can find a suitable $w$ in the isometry group of the model space $\A'$. We only need to show that one can find such a $w$ in the correct subgroup of this isometry group associated to the spherical Coxeter system $(\overline{W}',S')$. If two apartments of $(\Lambda',\cF')$ share no point then this is trivial, so assume they do share a point $\alpha$. Using the local structure around $\Lambda'_\alpha$, in particular Condition (B3), it follows that $w$ indeed can be found in the correct subgroup.

For Condition (A3), we start by proving a finite covering result for pairs of apartments.

\begin{lemma}
Let $\Sigma_1$ and $\Sigma_2$ be two apartments of $\Lambda'$. Let $\mathcal{C} := \mathcal{C}(\Sigma_1, \Sigma_2)$ be the set of apartments containing at infinity two chambers $c_i\in {\Sigma_i}_\infty$, $i=1,2$,  which are opposite in $\Lambda'_\infty$. Then $\mathcal{C}$ is a finite set of apartments such that each pair of points $(\alpha,\beta)\in \Sigma_1\times \Sigma_2$ is contained in one of these apartments in $\mathcal{C}$.
\end{lemma}
\proof
We first prove that if two points, one in each of the apartments $\Sigma_1$ and $\Sigma_2$, are contained in a common apartment, then they also lie in an apartment belonging  to $\mathcal{C}$. So suppose that $\Xi$ is an apartment containing points $\alpha_i\in \Sigma_i$, $i=1,2$. Let $S_{\alpha_1}$ be a sector of $\Xi$ based at $\alpha_1$ containing $\alpha_2$. The fact that $\Lambda'_{\alpha_1}$ is a spherical building together with~\cite[Exercise 4.90]{Abr-Bro:08} implies the existence of an $\alpha_1$ based sector $T_{\alpha_1}$ contained in $\Sigma_1$ whose germ is opposite the germ of $S_{\alpha_1}$ at $\alpha_1$. The set $T_{\alpha_1}\cup S_{\alpha_1}$ is contained in some apartment $\Xi'$ by Lemma~\ref{p12}. The sector $T_{\alpha_2}$ based at $\alpha_2$ parallel to $T_1$ now contains $\alpha_1$. 

Let $T'_{\alpha_2}$ be a sector of $\Sigma_2$ based at $\alpha_2$ and such that its germ is opposite the germ of $T_{\alpha_2}$. Again by Lemma~\ref{p12} we obtain a unique apartment $\Xi''$ containing $T_{\alpha_2}$ and $T'_{\alpha_2}$ and hence also the points $\alpha_1$ and $\alpha_2$. Since apartments of $(\Lambda',\cF')$ are in one to one correspondence with apartments in $\Lambda'_\infty$ the apartment $\Xi''$ is uniquely determined by its chambers at infinity $T_\infty$ and $T'_\infty$. By construction the apartment $\Xi''$ is contained in $\mathcal{C}$. So we conclude that if two points $\alpha_i\in \Sigma_i$, $i=1,2$ lie in a common apartment, then they also lie in a common apartment belonging to $\mathcal{C}$.

Now suppose that two points $\alpha_i\in \Sigma_i$ do not lie in one apartment. Let $K$ be the set of points in $\Sigma_1$ which do lie in an apartment with $\alpha_2$. Due to the above discussion and the fact that (A2) is already proven for $(\Lambda',\cF')$ we have that $K$ is a finite union of convex sets. Let $S$ be a sector of $\Sigma_1$. Lemma~\ref{p8} implies that there exists an apartment of $(\Lambda',\cF')$  containing a sector based at $\alpha_2$ and parallel to $S$. By the discussion in Section~\ref{section:subbie} it follows that this apartment and $\Sigma_1$ share at least a sector. So $K$ is not empty, but also not the entirety of $\Sigma_1$ as $\alpha_1 \notin K$. Because $K$ is a finite union of convex sets, one can find a point $\beta$ in $K$ so that not all the germs based at $\beta$ in $\Sigma_1$ lie in $K$ (a point of the boundary of $K$ in $\Sigma_1$). Let $R_\beta$ be a sector-facet with such a germ, and $U_\beta$ a sector based at $\beta$ containing $\alpha_2$ (possible because there exists an apartment containing both). Lemma~\ref{p8} yields that there exists an apartment $\Sigma$ containing the sector $U_\beta$ and the germ of $R_\beta$. The germ of $R_\beta$ now lies in $K$, contradicting the way we have chosen $R_\beta$. So we obtain that $K$ contains all points of $\Sigma_1$. The lemma is hereby proven.
\qed

\begin{cor}
Condition \emph{(A3)} is satisfied by $(\Lambda',\cF')$.\end{cor}
\proof
Directly from the above lemma.
\qed

As we have proven Conditions (A1)-(A5) and (TI), it follows that $(\Lambda',\cF)$ forms indeed an $\R$-building. Its building at infinity is the fixed structure of $G$ in $\Lambda_\infty$ by construction.

\subsection{Isometries acting on $\Lambda'$}
Let $G'$ be the group of isometries of $\Lambda$ centralizing $G$. We want to show that this group acts on $\Lambda'$. For this we first study the point sets of embedded $\R$-buildings with building at infinity $\Lambda'$. Let $S_\infty$ be a maximal stabilized simplex of $\Lambda_\infty$.

\begin{lemma}
Let $K$ be the union of the point sets of each $\R$-building embedded in fixed point set of $G$ in $\overline{\Lambda}$ with $\Lambda'_\infty$ as building at infinity. Then $K$ is isometric to $\Lambda' \times K'$, where $K'$ is a convex subset of $T(S_\infty)$.
\end{lemma}
\proof
This follows from Section 2.3.3 of~\cite{Kle-Lee:97}.
\qed

The group $G'$ stabilizes $\Lambda'_\infty$ and the fixed point set of $G$ in $\overline{\Lambda}$, so it stabilizes $K$. By dividing out the factor $K'$, we obtain an action of $G'$ on $\Lambda'$.

By repeating the above arguments twice (see the introduction of Section~\ref{section:proof}), one has proven the main result.

\begin{remark} \rm
If we are in the case where $F$ is bounded as described in Remark~\ref{rem:applibound}, then $G'$ stabilizes the set of points $\Lambda'$ directly.
\end{remark}

\section{Tits diagrams with minimal angle strictly greater than $\pi/3$}\label{subsection:angle}
In the previous sections we constructed under certain conditions an $\R$-building with $\Lambda'_\infty$ as building at infinity. One of these conditions involves the minimal angle of the Tits diagram. In this section we list some Tits diagrams with  minimal angle strictly greater than $\pi/3$.

We only need to investigate the case where the Tits diagram $(M,\Gamma,A)$ is of relative rank one and has trivial $\Gamma$ (see Section~\ref{section:angle}). The most interesting diagrams are those arising from the classification of algebraic semi-simple groups in~\cite{Tit:66}. 

Note that when the diagram is disconnected, one only needs to consider the connected component with the isotropic orbit.

\paragraph{Minimal angle $\pi$.}

If the diagram only consists of the isotropic orbit, it is clear that the possible angles are 0 and $\pi$. This occurs in the quasi-split case mentioned in the introduction.

\paragraph{Minimal angle $\pi/2$.}
Examples here are the diagrams of type $\mathsf{B}_n$ or $\mathsf{D}_n$ where the first node is the isotropic orbit, so: 
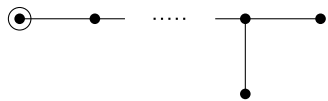
\begin{center}
\begin{tikzpicture}
\fill (-2,0) circle (2pt);
\fill (-1,0) circle (2pt);
\fill (1,0) circle (2pt);
\fill (2,0) circle (2pt);

\draw (-2,0) circle (.15);
\draw[dotted, thick] (-.22,0) -- (.22,0); 

\draw (-2,0) -- (-.6,0); 
\draw (.6,0) -- (1,0); 

\draw (1,.05) -- (2,.05); 
\draw (1,-.05) -- (2,-.05); 

\end{tikzpicture}
\end{center}
and
\begin{center}
\begin{tikzpicture}
\fill (-2,0) circle (2pt);
\fill (-1,0) circle (2pt);
\fill (1,0) circle (2pt);
\fill (1,-1) circle (2pt);

\fill (2,0) circle (2pt);

\draw (-2,0) circle (.15);
\draw[dotted, thick] (-.22,0) -- (.22,0); 
\draw (-2,0) -- (-.6,0); 
\draw (.6,0) -- (2,0); 

\draw (1,0) -- (1,-1);

\end{tikzpicture}
\end{center}

Using the notations of~\cite[Ch. VI, \S 4]{Bou:68} one checks that in both cases the vertices of the isotropic orbit correspond to the vectors in $\R^n$ of the form $(0,\dots,0, \pm 1,0,\dots,0)$. From this it follows that the possible angles are $0$, $\pi/2$ and $\pi$.

Another example is the diagram of type $\mathsf{A}_3$ where the middle node is encircled, see Section~\ref{section:diag}.

\paragraph{Minimal angle $\arccos(1/3)$.}
These angles occur for the following Tits diagrams of type $\mathsf{A}_5$ and $\mathsf{E}_7$:
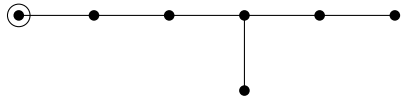
\begin{center}
\begin{tikzpicture}
\fill (-2,0) circle (2pt);
\fill (-1,0) circle (2pt);
\fill (0,0) circle (2pt);
\fill (1,0) circle (2pt);
\fill (2,0) circle (2pt);

\draw (-2,0) -- (2,0); 

\draw (0,0) circle (.15);

\end{tikzpicture}, 

\end{center}
and
\begin{center}
\begin{tikzpicture}
\fill (-2,0) circle (2pt);
\fill (-1,0) circle (2pt);
\fill (0,0) circle (2pt);
\fill (1,0) circle (2pt);
\fill (2,0) circle (2pt);
\fill (3,0) circle (2pt);
\fill (1,-1) circle (2pt);
\draw (1,0) -- (1,-1); 

\draw (-2,0) -- (3,0); 

\draw (-2,0) circle (.15);

\end{tikzpicture}
\end{center}

These cases also appear in the theory of line systems with angle $\alpha$ such that $\cos(\alpha)= \pm 1/3$, see examples $C_{10}$ and $C_{28}$ in~\cite{Shu-Yan:80}. Alternatively one can calculate the vertices and verify the angles directly. For the $\mathsf{A}_5$ diagram the vertices correspond to vectors obtained from permutations of $(1,1,1,-1,-1,-1)$, and for the $\mathsf{E}_7$ diagram they correspond to the vectors $(2,0,0,0,0,0,1,-1)$, $(-2,0,0,0,0,0,1,-1)$, $(1,1,1,1,1,1,0,0)$ and $(-1,-1,1,1,1,1,0,0)$ up to permutation of the first six coordinates and taking the negative vector (again using the notations of~\cite[Ch. VI, \S 4]{Bou:68}).

Other examples are the following Tits diagrams of type $\mathsf{B}_3$ and $\mathsf{A}_5$:
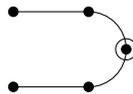
\begin{center}
\begin{tikzpicture}
\fill (2,0) circle (2pt);
\fill (1,0) circle (2pt);
\fill (3,0) circle (2pt);
\draw (2,.035) -- (3,.035); 
\draw (2,-.035) -- (3,-.035); 
\draw (3,0) circle (.15);
\draw (1,0) -- (2,0); 
\end{tikzpicture}
\end{center}
and 
\begin{center}\begin{tikzpicture}
\fill (0,0) circle (2pt);
\fill (1,1) circle (2pt);
\fill (1,0) circle (2pt);
\fill (0,1) circle (2pt);
\fill (1.5,.5) circle (2pt);

\draw (1.5,.5) circle (4pt);
\draw (0,0) -- (1,0); 
\draw (0,1) -- (1,1); 

\draw (1,0) arc (-90:90:.5);

\end{tikzpicture}\end{center}
Note that the first diagram is the folded version of the second one. That the minimal angle is indeed $\arccos(1/3)$ is easily verified by direct calculus.

\section{Existence of $\R$-buildings corresponding to certain forms}\label{exis}

In this section we discuss an application of the main result to the existence of $\R$-buildings for Moufang polygons of exceptional type (see Section~\ref{section:bruhattits}). It has been conjectured by Jacques Tits in~\cite[p. 173]{Tit:86} that, given a Moufang polygon (of exceptional type) defined over a field with a complete valuation which is preserved by an involution or Tits-endomorphism $\sigma$ where applicable, there exists an $\R$-building with this Moufang polygon as building at infinity.
This existence has been already proven algebraically for all Moufang polygons of exceptional type except for the quadrangle of type $\mathsf{E}_8$ (see~\cite{Wei:09} and the recent~\cite{Wei:*}). We now show how our results can be applied to reprove this in a geometric way, except for the quadrangle of type $\mathsf{E}_7$.

Consider a spherical building of type $\mathsf{F}_4$ or $\mathsf{E}_i$ ($i= 6,7,8$), defined over a field with a complete valuation. One can realize such a building as the spherical building $\Lambda_\infty$ at infinity of an $\R$-building $(\Lambda,\cF)$. Moufang polygons of exceptional type now arise as fixed structures $\Lambda'_\infty$ of finite groups of isomorphisms of $\Lambda_\infty$. Due to completeness of the valuation this group action extends to the $\R$-building.

If the minimal angle of the Tits diagram is strictly greater than $\pi/3$ one can now apply our main result and obtain an $\R$-building $(\Lambda',\cF')$ with the generalized polygon as building at infinity. In particular this shows the existence of $\R$-buildings where the Moufang polygon of exceptional type has the following Tits diagram.

\begin{itemize} 
\item Generalized triangles:

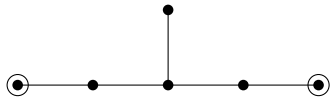
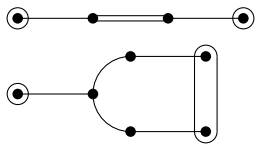
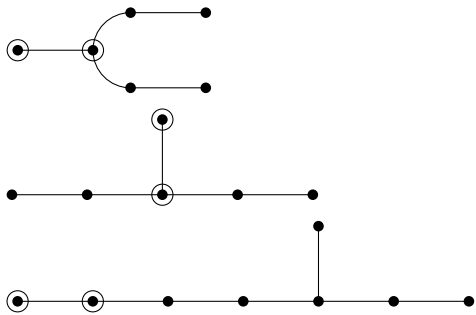
\begin{tikzpicture}
\foreach \x in {1,...,5}{
\fill (\x,0) circle (2pt);}
\draw (1,0) circle (4pt);
\draw (5,0) circle (4pt);

\fill (3,1) circle (2pt);

\draw (1,0) -- (5,0); 
\draw (3,0) -- (3,1); 

\end{tikzpicture}

\item Generalized quadrangles:

\begin{tikzpicture}
\foreach \x in {1,2,3,4}{
\fill (\x,0) circle (2pt);}
\draw (1,0) circle (4pt);
\draw (4,0) circle (4pt);

\draw (1,0) -- (2,0); 
\draw (4,0) -- (3,0); 
\draw (2,.035) -- (3,.035); 
\draw (2,-.035) -- (3,-.035); 

\end{tikzpicture}

\begin{tikzpicture}
\fill (0,0) circle (2pt);
\fill (-1,1) circle (2pt);
\fill (-1,0) circle (2pt);
\fill (0,1) circle (2pt);
\fill (-1.5,.5) circle (2pt);
\fill (-2.5,.5) circle (2pt);

\draw (-2.5,.5) circle (4pt);

\draw (0,0) -- (-1,0); 
\draw (0,1) -- (-1,1); 
\draw (-1.5,.5) -- (-2.5,.5);

\draw (-1,1) arc (90:270:.5);
\draw (-.15,0) arc (-180:0:.15);
\draw (-.15,1) arc (180:0:.15);
\draw (-.15,0) -- (-.15,1); 
\draw (.15,0) -- (.15,1); 

\end{tikzpicture}

\item Generalized hexagons:

\begin{tikzpicture}
\fill (0,0) circle (2pt);
\fill (-1,1) circle (2pt);
\fill (-1,0) circle (2pt);
\fill (0,1) circle (2pt);
\fill (-1.5,.5) circle (2pt);
\fill (-2.5,.5) circle (2pt);

\draw (-2.5,.5) circle (4pt);
\draw (-1.5,.5) circle (4pt);
\draw (0,0) -- (-1,0); 
\draw (0,1) -- (-1,1); 
\draw (-1.5,.5) -- (-2.5,.5);

\draw (-1,1) arc (90:270:.5);

\end{tikzpicture}

\begin{tikzpicture}
\foreach \x in {1,...,5}{
\fill (\x,0) circle (2pt);}
\draw (3,0) circle (4pt);
\draw (3,1) circle (4pt);

\fill (3,1) circle (2pt);

\draw (1,0) -- (5,0); 
\draw (3,0) -- (3,1); 

\end{tikzpicture}

\begin{tikzpicture}
\foreach \x in {0,1,...,6}{
\fill (\x,0) circle (2pt);}
\draw (0,0) circle (4pt);
\draw (1,0) circle (4pt);

\fill (4,1) circle (2pt);

\draw (0,0) -- (6,0); 
\draw (4,0) -- (4,1); 

\end{tikzpicture}

\end{itemize}

The Moufang quadrangles of exceptional type $\mathsf{E}_7$ and $\mathsf{E}_8$ have minimal angle $\pi/3$, so our present method does not yield existence for these. But these, and also some rank one buildings, will be handled in Part II of this paper.


\end{document}